\documentclass[10pt]{amsart}
\usepackage{amsmath,amsfonts,amssymb}
\usepackage{enumerate}

\newtheorem{lemma}{Lemma}[section]

\newtheorem{theorem}[lemma]{Theorem}
\newtheorem{corollary}[lemma]{Corollary}
\theoremstyle{definition}

\newtheorem{question}[lemma]{Question}
\theoremstyle{remark}

\numberwithin{equation}{section} \numberwithin{table}{section}

\newcommand{\rank}{\mathrm{rank}\,}

\begin{document}
\title[Rank one matrices and the finiteness property]{Rank one matrices do not contribute to the failure of the finiteness property}
\author{Ian D. Morris}
\address{Dipartimento di Matematica, Universit\'a degli Studi di Roma Tor Vergata, Via della Ricerca Scientifica, Rome 00133, Italy.}
\email{ian.morris.ergodic@gmail.com}

\begin{abstract}
The joint spectral radius of a bounded set of $d \times d$ real or complex matrices is defined to be the maximum exponential rate of growth of products of matrices drawn from that set. A set of matrices is said to satisfy the finiteness property if this maximum rate of growth occurs along a periodic infinite sequence. In this note we give some sufficient conditions for a finite set of matrices to satisfy the finiteness property in terms of its rank one elements. We show in particular that if a finite set of matrices does not satisfy the finiteness property, then the subset consisting of all matrices of rank at least two is nonempty, does not satisfy the finiteness property, and has the same joint spectral radius as the original set. We also obtain an exact formula for the joint spectral radii of sets of matrices which contain at most one element not of rank one, generalising a recent result of X. Dai.
\end{abstract}

\maketitle

\section{Introduction and statement of results}

Given a bounded set $\mathsf{A}$ of $d \times d$ matrices over $\mathbb{R}$ or $\mathbb{C}$, the joint spectral radius of $\mathsf{A}$ is defined to be the quantity
\[\varrho(\mathsf{A})=\lim_{n \to \infty} \sup\left\{\left\|A_{i_n} \cdots A_{i_1}\right\|^{\frac{1}{n}} \colon A_i \in \mathsf{A}\right\},\]
a definition introduced by G.-C. Rota and G. Strang in 1960 (\cite{RS}, subsequently reprinted in \cite{Rotacoll}). This limit always exists and is independent of the norm used (for a proof see e.g. \cite{Jungers}). The joint spectral radius arises naturally in a range of mathematical contexts including control and optimisation, wavelet regularity, coding theory, and combinatorics; for references see \cite{Chang,HMST,Jungers,PJB}. As such the properties of the joint spectral radius are the subject of ongoing research investigation \cite{AAJPR,Chang, GuPr,HMST,Kozb2}.

Given a bounded set of matrices $\mathsf{A}$, let us define $\mathsf{A}_n:=\{A_{i_1}\cdots A_{i_n} \colon A_j \in \mathsf{A}\}$ for every integer $n \geq 1$, and let $\rho(A)$ denote the ordinary spectral radius of a square matrix $A$. The joint spectral radius satisfies the identity
\[\varrho(\mathsf{A}) = \limsup_{n \to \infty} \sup_{A \in \mathsf{A}_n}\rho(A)^{1/n}\]
established by M. A. Berger and Y. Wang in \cite{BW}, and I. Daubechies and J. Lagarias noted in \cite{DL} that the inequalities
\[ \sup_{A \in \mathsf{A}_n}\rho(A)^{1/n} \leq \varrho(\mathsf{A}) \leq  \sup_{A \in \mathsf{A}_n}\|A\|^{1/n}\]
apply for every $n \geq 1$, where $\| \cdot\|$ is any operator norm on the space of $d \times d$ matrices. The joint spectral radius therefore admits the characterisation
\begin{align}\label{basic}\varrho(\mathsf{A})&=\lim_{n \to \infty} \sup_{A \in \mathsf{A}_n} \|A\|^{1/n} = \inf_{n \geq 1} \sup_{A \in \mathsf{A}_n} \|A\|^{1/n}\\
&=\limsup_{n \to \infty} \sup_{A \in \mathsf{A}_n}\rho(A)^{1/n} = \sup_{n \geq 1}\sup_{A \in \mathsf{A}_n} \rho(A)^{1/n}.\nonumber\end{align}
for every operator norm $\|\cdot\|$. Motivated by these identities, J. Lagarias and Y. Wang asked in \cite{LW} whether every finite set of $d \times d$ matrices $\mathsf{A}$ has the property that there exists an integer $N$ depending on $\mathsf{A}$ such that $\varrho(\mathsf{A})=\max\{\rho(A)^{1/N} \colon A \in \mathsf{A}_N\}$. If such an integer exists then $\mathsf{A}$ is said to have the \emph{finiteness property}. Various authors have since shown that there exist pairs of invertible $2 \times 2$ real matrices which do not have the finiteness property \cite{BTV,BM,HMST,Koz3}. 
Conversely, some sufficient conditions for the finiteness property to hold have been given in \cite{Daiquack,DaiKoz,Gu,Gu2,LW}, and the finiteness property for rational matrices has been studied in detail in \cite{CGSC,BJ}. In this note we investigate the relationship between the finiteness property of a set $\mathsf{A}$ and the presence of rank one matrices in that set. We will show in particular that the known examples of pairs of matrices which lack the finiteness property are in some sense the simplest possible: if a finite set of $d \times d$ matrices fails to have the finiteness property, then it must be the case that at least two of its elements have rank greater than or equal to two.

The main result of this note is the following theorem:
\begin{theorem}\label{bum}
Let $\mathsf{A}=\{A_1,\ldots,A_\ell\}$ be a finite set of $d \times d$ matrices over $\mathbb{C}$, and let $\mathsf{R} \subseteq \mathsf{A}$ be the set of all elements of $\mathsf{A}$ which have rank greater than or equal to two. If  $\varrho(\mathsf{R})<\varrho(\mathsf{A})$ or $\mathsf{R}$ is empty, then there exists a finite sequence $(A_{i_j})_{j=1}^n$ of elements of $\mathsf{A}$, in which each rank-one element of $\mathsf{A}$ appears at most once, such that $\varrho(\mathsf{A})=\rho(A_{i_1}\cdots A_{i_n})^{1/n}$.
\end{theorem}
The following corollary motivates the title of this article. (Note that a singleton set trivially satisfies the finiteness property by Gelfand's formula.)
\begin{corollary}\label{boo}

Let $\mathsf{A}$ be a finite set of $d \times d$ matrices over $\mathbb{C}$,  and let $\mathsf{R} \subseteq \mathsf{A}$ be the set of all elements of $\mathsf{A}$ which have rank greater than or equal to two. If $\mathsf{A}$ does not satisfy the finiteness property then $\varrho(\mathsf{A})=\varrho(\mathsf{R})$, $\mathsf{R}$ contains at least two elements, and $\mathsf{R}$ does not satisfy the finiteness property.
\end{corollary}
 Note that the converse of Corollary \ref{boo} is false: if $\mathsf{B}:=\{B_1,B_2\}$ is a pair of invertible $2 \times 2$ real matrices which fails to satisfy the finiteness property, then the set of $3 \times 3$ matrices defined by
\[\mathsf{A}:=\left\{\left(\begin{array}{cc}B_1 & 0 \\ 0 &0\end{array}\right), \left(\begin{array}{cc}B_2 & 0 \\ 0 &0\end{array}\right), \left(\begin{array}{cc}0 & 0 \\ 0 &\varrho(\mathsf{B})\end{array}\right)\right\}\]
clearly satisfies the finiteness property even though $\varrho(\mathsf{A})=\varrho(\mathsf{R})$ and $\mathsf{R}$ does not  satisfy the finiteness property.

In the special case where $\mathsf{R}$ contains at most one element, Theorem \ref{bum} reduces to the following:
\begin{corollary}\label{waxman}
Let $\mathsf{A}=\{A_1,\ldots,A_\ell\}$ be a finite set of $d \times d$ matrices over $\mathbb{C}$, and suppose that $\mathsf{A}$ contains at most one element with rank greater than or equal to two. Then there exists a finite sequence $(A_{i_j})_{j=1}^n$ of elements of $\mathsf{A}$, in which each rank-one element of $\mathsf{A}$ appears at most once, such that $\varrho(\mathsf{A})=\rho(A_{i_1}\cdots A_{i_n})^{1/n}$.
\end{corollary}
In the special case $\ell=2$ this further reduces to the following result, which strengthens a recent theorem of X. Dai \cite{Daiquack}.
\begin{corollary}
If $\mathsf{A}=\{A,B\}$ where $A$ has rank one, then either $\varrho(\mathsf{A})=\rho(B)$ or $\varrho(\mathsf{A})=\rho(AB^n)^{1/(n+1)}$ for some integer $n \geq 0$. 
\end{corollary}
It should be noted that the stipulation in Theorem \ref{bum} that $\mathsf{A}$ should be finite rather than compact is deliberate, since the proof of Theorem \ref{bum} makes use of the pigeonhole principle applied to the elements of $\mathsf{A}$. In particular it remains unclear whether an analogue of Corollary \ref{boo} holds when $\mathsf{A}$ is compact and infinite. Motivated by this observation we ask the following question:
\begin{question}
Does there exist a compact infinite set of rank one matrices which does not have the finiteness property?
\end{question}

\section{Proof of Theorem \ref{bum}}
If $\varrho(\mathsf{A})=0$ then the theorem holds trivially, since it follows from \eqref{basic} that $\rho(A)=0=\varrho(\mathsf{A})$ for every $A \in \mathsf{A}$. We may therefore freely assume that $\varrho(\mathsf{A})$ is positive. By multiplying every element of $\mathsf{A}$ by the nonzero real scalar $\varrho(\mathsf{A})^{-1}$, we may further reduce to the case in which $\varrho(\mathsf{A})=1$. For the rest of the proof we shall assume without loss of generality that $\varrho(\mathsf{A})=1$.

Given a finite set of $d \times d$ complex matrices $\mathsf{B}=\{B_1,\ldots,B_k\}$ with nonzero joint spectral radius, we shall say that a norm $\|\cdot\|$ on $\mathbb{C}^d$ is an \emph{extremal norm} for $\mathsf{B}$ if the induced operator norm on the space of $d \times d$ matrices satisfies $\|B\| \leq \varrho(\mathsf{B})$ for all $B \in \mathsf{B}$. We will begin the proof by showing that Theorem \ref{bum} holds subject to the additional assumption that $\mathsf{A}$ admits an extremal norm.

We will prove the theorem in the case where $\mathsf{R}$ is nonempty: the case where this set is empty may be dealt with by a simple modification of the following argument. Let us then write $\mathsf{A} = \{A_1,A_2,\ldots,A_\ell\}$ with $\mathsf{R}=\{A_1,\ldots,A_r\} \subset \mathsf{A}$. Using \eqref{basic} we note that there exists $N \geq 1$ such that
\[\sup_{k \geq N} \sup \{\|B\| \colon B \in \mathsf{R}_k\}<1.\]
Let $n \geq N(1+\ell-r)$ and, using \eqref{basic} together with the fact that $\|\cdot\|$ is extremal for $\mathsf{A}$, choose $i_1,\ldots,i_n$ such that $\|A_{i_1} \cdots A_{i_n}\|=1$. If $m, k \geq 1$ are integers such that $i_{m},i_{m+1},\ldots,i_{m+k-1}$ all belong to the set $\{1,\ldots, r\}$, then we necessarily have $k < N$ since
\begin{align*}1 = \left\|A_{i_1} \cdots A_{i_n}\right\| &\leq \left\|A_{i_1} \cdots A_{i_{m-1}}\right\| . \left\|A_{i_m} \cdots A_{i_{m+k-1}}\right\| . \left\|A_{i_{m+k}} \cdots A_{i_n}\right\| \\
&\leq \left\|A_{i_m} \cdots A_{i_{m+k-1}}\right\| \leq \sup_{B \in \mathsf{R}_k}\|B\|\end{align*}
and this last term is greater than or equal to $1$ only when $k<N$. In particular, each block of $N$ successive symbols $i_j,\ldots,i_{j+N-1}$ must include at least one entry in the range $\{r+1,\ldots,\ell\}$.  Since $n \geq N(1+\ell-r)$ the sequence $i_1,\ldots,i_n$ includes at least $1+\ell-r$ pairwise non-overlapping blocks of length $N$, and it follows by the pigeonhole principle that there exist integers $k_1,k_2$ such that $1 \leq k_1 < k_2 \leq n$ and $A_{i_{k_1}} = A_{i_{k_2}} \notin \mathsf{R}$. Choose integers $k_1,k_2$ with these properties such that the difference $k_2-k_1$ is minimised. If $A_{i_{\ell_1}}=A_{i_{\ell_2}} \notin\mathsf{R}$ and $k_1 \leq \ell_1 <\ell_2 \leq k_2$ then necessarily $\ell_1=k_1$ and $\ell_2=k_2$ by minimality, and it follows that those elements of the sequence $(A_{i_j})_{j=k_1}^{k_2-1}$ which belong to $\mathsf{A} \setminus \mathsf{R}$ are necessarily all distinct.

We claim that $\rho(A_{i_{k_1}} \cdots A_{i_{k_2-1}})=1$, which proves the theorem in the case of the matrix set $\mathsf{A}$. Since
\begin{align*}
1=\|A_{i_1}\cdots A_{i_n}\| &\leq \|A_{i_1} \cdots A_{i_{k_1-1}}\|. \|A_{i_{k_1}} \cdots A_{i_{k_2}}\| . \|A_{i_{1+k_2}} \cdots A_n\|\\ &\leq \left\|A_{i_{k_1}} \cdots A_{i_{k_2}}\right\| \leq 1\end{align*}
we have $\left\|A_{i_{k_1}} \cdots A_{i_{k_2}}\right\|=1$. Choose a vector $v$ such that $\left\|A_{i_{k_1}} \cdots A_{i_{k_2}}v\right\|=1$ and $\|v\|=1$ . Since
\[1=\left\|A_{i_{k_1}}\cdots A_{i_{k_2}}v \right\| \leq \left\|A_{i_{k_1}} \cdots A_{i_{k_2-1}}\right\|.\left\|A_{i_{k_2}}v\right\| \leq \left\|A_{i_{k_2}}v\right\| \leq 1\] 
we have $\left\|A_{i_{k_2}}v\right\|=1$. The vector $A_{i_{k_2}}v$ thus has unit length and belongs to the image of
 the matrix $A_{i_{k_2}}$. The product $A_{i_{k_1}} \cdots A_{i_{k_2-1}}$ maps this vector to a vector which also has unit length and also belongs to the image of $A_{i_{k_2}}=A_{i_{k_1}}$. Since $A_{i_{k_2}}$ has rank one its image is one-dimensional, and it follows that $v$ is an eigenvector of $A_{i_{k_1}} \cdots A_{i_{k_2-1}}$ whose corresponding eigenvalue has unit modulus. We therefore have $\rho(A_{i_{k_1}} \cdots A_{i_{k_2-1}})=1$ as claimed, completing the proof of the theorem in this case.

It remains to show that this result implies the validity of Theorem \ref{bum} in general. To achieve this we use a reduction argument which is by now somewhat standard. We require the following lemma, the proof of which may be found in e.g. \cite{E}.
\begin{lemma}\label{els}
Suppose that the set of $d \times d$ matrices $\mathsf{A}=\{A_1,\ldots,A_\ell\}$ has nonzero joint spectral radius and does not admit an extremal norm. Then there exist positive integers $d_1,d_2$ such that $d_1+d_2=d$, a family of $d_1 \times d_1$ matrices $\mathsf{B}=\{B_1,\ldots,B_\ell\}$, a family of $d_1 \times d_2$ matrices $\mathsf{C}=\{C_1,\ldots,C_\ell\}$, a family of $d_2 \times d_2$ matrices $\mathsf{D}=\{D_1,\ldots,D_\ell\}$ and an invertible $d \times d$ matrix $U$ such that
\begin{equation}\label{dux}U^{-1}A_iU = \left(\begin{array}{cc}B_i&C_i\\0&D_i\end{array}\right)\end{equation}
for $i=1,\ldots,\ell$.
\end{lemma}
Note that the matrices $B_i$ are not necessarily pairwise distinct, and similarly for the families of matrices $C_i$ and $D_i$.

We will prove the full statement of Theorem \ref{bum} by induction on the dimension $d$ of the matrices which comprise the set $\mathsf{A}$. Since $1 \times 1$ matrix multiplication is commutative, it is not difficult to see that the joint spectral radius of a finite set of $1 \times 1$ matrices is just the maximum of the spectral radii of the individual matrices, and in particular Theorem \ref{bum} holds trivially in this case. Let us now suppose that Theorem \ref{bum} has been shown to hold for all finite sets of square matrices of dimension at most $d$, and use this fact to establish that it is also valid for sets of square matrices of dimension $d+1$.

Let $\mathsf{A}=\{A_1,\ldots,A_\ell\}$ be a finite set of $(d+1)\times(d+1)$ matrices which meets the hypotheses of Theorem \ref{bum}. If $\mathsf{A}$ admits an extremal norm then Theorem \ref{bum} is valid for $\mathsf{A}$ by the argument given previously. Otherwise, let $d_1,d_2,\mathsf{B}, \mathsf{C}, \mathsf{D}$ and $U$ be as in Lemma \ref{els}. An easy consequence of \eqref{dux} is that $\rho(A_i)=\max\{\rho(B_i),\rho(D_i)\}$ for each $i=1,\ldots,\ell$, and in the same manner one may show that $\rho(A_{i_1} \cdots A_{i_n})=\max\{\rho(B_{i_1}\cdots B_{i_n}),\rho(D_{i_1}\cdots D_{i_n})\}$ for any finite sequence $i_1,\ldots,i_n$ of elements of $\{1,\ldots,\ell\}$. Using \eqref{basic} we deduce that $\varrho(\mathsf{A})=\max\{\varrho(\mathsf{B}),\varrho(\mathsf{D})\}$. We suppose that $\varrho(\mathsf{A})=\varrho(\mathsf{B})$, the proof in the complementary case $\varrho(\mathsf{A})=\varrho(\mathsf{D})$ being similar. 

Let $\mathsf{R}_{\mathsf{B}}$ be the set of all elements of $\mathsf{B}$ which have rank greater than or equal to two. If $\mathsf{R}_{\mathsf{B}}$ is empty, then the hypotheses of Theorem \ref{bum} apply to $\mathsf{B}$. If $\mathsf{R}_{\mathsf{B}}$ is not empty, then since each matrix $B_i$ is obtained from the corresponding matrix $A_i$ by composition with an invertible matrix and a projection we necessarily have $\rank B_i \leq \rank A_i$. In particular if $B_i$ has rank at least two, then $A_i$ does also. We deduce that
\[\max\left\{\rho(B_{i_1}\cdots B_{i_n})^{1/n} \colon B_j \in \mathsf{R}_{\mathsf{B}}\right\} \leq \max\left\{\rho(A_{i_1}\cdots A_{i_n})^{1/n} \colon A_j \in \mathsf{R}\right\}\]
for every $n \geq 1$, and using \eqref{basic} we may derive the inequality $\varrho(\mathsf{R}_{\mathsf{B}}) \leq \varrho(\mathsf{R})<\varrho(\mathsf{A})=\varrho(\mathsf{B})$. The hypotheses of Theorem \ref{bum} therefore also apply to $\mathsf{B}$ in the case where $\mathsf{R}_{\mathsf{B}}$ is not empty. Since the matrices comprising $\mathsf{B}$ have dimension $d_1= d+1 - d_2 \leq d$ it follows by the induction hypothesis that Theorem \ref{bum} is valid for $\mathsf{B}$, and hence there exists a finite sequence $({i_j})_{j=1}^n$ of integers in the range $1,\ldots,\ell$ such that $\varrho(\mathsf{B})=\rho(B_{i_1}\cdots B_{i_n})^{1/n}$ and such that no rank-one matrix $B_i$ occurs twice in the sequence $(B_{i_j})_{j=1}^n$.
Since $\varrho(\mathsf{A})=\varrho(\mathsf{B})$ and $\rank B_i \leq \rank A_i$ for each $i$, we deduce that $\varrho(\mathsf{A})=\rho(A_{i_1}\cdots A_{i_n})^{1/n}$ and no rank-one matrix $A_i$ occurs twice in the sequence $(A_{i_j})_{j=1}^n$. We conclude that Theorem \ref{bum} is valid for $\mathsf{A}$, which completes the proof of the induction step. This completes the proof of Theorem \ref{bum} in full generality.

\section{Acknowledgment}

This research was conducted as part of the ERC grant MALADY (246953).

\bibliographystyle{siam}
\bibliography{Rcon}

\end{document}